\documentclass{ITaSconf}
\usepackage[utf8x]{inputenc}
\usepackage[english]{babel}
\usepackage{url}

\usepackage{etex}
\usepackage{amssymb}
\usepackage{euscript}
\usepackage{upgreek}
\usepackage{array}
\usepackage{theorem}
\usepackage{graphicx}
\usepackage{subfig}
\usepackage{comment}
\usepackage{cite}
\usepackage{gensymb} 
\usepackage{amsmath}

\newtheorem{Def}{Definition}

\newcommand{\bx}{\mathbf{x}}

\newcommand{\by}{\mathbf{y}}
\newcommand{\bb}{\mathbf{b}}

\newcommand{\bY}{\mathbf{Y}}
\newcommand{\bA}{\mathbf{A}}

\newcommand{\bQ}{\mathbf{Q}}

\newcommand{\bbR}{\mathbb{R}}

\newcommand{\T}{^{\text{\tiny\sffamily\upshape\mdseries T}}}
\newcommand{\argmin}{\mathop{\arg \min}\limits}
\newcommand{\argmax}{\mathop{\arg \max}\limits}
\newcommand\undermat[2]{%
  \makebox[0pt][l]{$\smash{\underbrace{\phantom{%
    \begin{matrix}#2\end{matrix}}}_{\text{$#1$}}}$}#2}

\title{Quadratic Programming Approach to Fit Protein Complexes into Electron Density Maps}
\author{
  Roman Pogodin\\
  \begin{affiliation}
  Skolkovo Institute of Science and Technology, Nobel St., 3, Moscow, 143026, Russia\\
    Moscow Institute of Physics and Technology, Institutskiy Lane 9, Dolgoprudny, Moscow, 141700, Russia\\ 
  \end{affiliation}\\
  \email{pogodin@phystech.edu}
  \and
  Alexander Katrutsa \\
  \begin{affiliation}
    Moscow Institute of Physics and Technology, Institutskiy Lane 9, Dolgoprudny, Moscow, 141700, Russia\\ Skolkovo Institute of Science and Technology, Nobel St., 3, Moscow, 143026, Russia\\
  \end{affiliation}\\
  \email{aleksandr.katrutsa@phystech.edu}
  \and
  Sergei Grudinin \\
  \begin{affiliation}
    University of Grenoble Alpes, LJK, F-38000 Grenoble, France\\
    CNRS, LJK, F-38000 Grenoble, France\\ Inria, F-38000 Grenoble, France
  \end{affiliation}\\
  \email{sergei.grudinin@inria.fr}
}

\begin{document}
\maketitle
\begin{abstract}
The paper investigates the problem of fitting protein complexes into electron density maps. They are represented by high-resolution cryoEM density maps converted into overlapping matrices and partly show a structure of a complex. The general purpose is to define positions of all proteins inside it. This problem is known to be NP-hard, since it lays in the field of combinatorial optimization over a set of discrete states of the complex. We introduce quadratic programming approaches to the problem. To find an approximate solution, we convert a density map into an overlapping matrix, which is generally indefinite. Since the matrix is indefinite, the optimization problem for the corresponding quadratic form is non-convex. 
To treat non-convexity of the optimization problem, we use different convex relaxations to find which set of proteins minimizes the quadratic form best.

\end{abstract}
\section{Introduction}
The problem of proteins fitting into cryoEM density maps of protein complexes remains important for biophysical studies of cell processes. Some examples of its importance can be found in \cite{Lawson2011}, which presents EMDataBank, an online database for electron microscopy. 

Two approaches to the problem are noticeable. The first one \cite{pandurangan} uses a genetic algorithm that discovers and then recombines good solutions which fit the density map. This parallel approach increases efficiency, but the accuracy decreases with the number of components inside the complex and the map's resolution. The second solution \cite{lasker} uses a cryoEM map directly and uses a set of predefined possible positions. It divides a set of fitting variables into uncoupled subsets, solves combinatorial optimization problems independently and finally gather all the solutions into the global minimum. This approach reduces the size of the problem from exponential in the number of all components to exponential in the number of components of the largest subset. However, it is sensitive to the accuracy of the component models and clustering into the subsets. 
Methods investigated in the current paper use a set of predefined positions like in the last paper. However, the minimization problem over a binary set is relaxed into a problem over a continuous set. The problem is then solved with continuous optimization methods, which are more efficient that discrete ones. The solution is then rounded to a binary one.

For continuous optimization methods, quadratic programming approaches and a stochastic algorithm are used in this paper. General overview of quadratic programming, main ideas and results of convex optimization are presented in \cite{boyd2004convex}.
For our purposes, the first idea is convex relaxation. The paper \cite{d2003relaxations} is a fulfilling review of a basic methods semidefinite (SDP). The second idea is sequential quadratic programming which finds a local minimum of a problem \cite{nocedal2006numerical}. The stochastic algorithm is called Simulated Annealing which also finds a local minimum of the problem \cite{Ingber199632}.

Efficiency of the method is tested in two principle ways. Firstly, it is testing against artificial cryoEM maps, which are based on a known structure of a protein. This method allows one to test efficiency at different map resolutions and was used in \cite{lasker, DiMaio2009181, Zhang2010, Pandurangan2012520, pandurangan, Topf2008295}. 
The second way is to use experimental density maps of protein complexes which structure we know, e.g. in \cite{lasker, Pandurangan2012520, pandurangan, Topf2008295}.

\section{Problem statement}
A protein complex consists of $m$ proteins and has $N$ computed spatial positions for each protein, which are different for different proteins. To predict the structure of the protein, we look for positions of a given set of proteins which fit into the density map best. Introduce correlation parameters for proteins' density maps and state the prediction problem formally.

\begin{Def}
Overlapping of two proteins' positions is overlapping of corresponding electron density maps.
\end{Def}
To measure overlapping between two density maps, we use the cross-correlation function \cite{Vasishtan2011333}, CCF:
\begin{equation}
\mathrm{CCF}=\underset{i}{\sum}\,\rho^1_i\rho_i^2,
\label{eq::ccf}
\end{equation}
where $\rho_i^1$ and $\rho_i^2$ are densities of $i$-th element of two maps. Laplacian-filtered CCF (LAP) allows one to compare matching of maps' edges rather than the whole volumes \cite{Vasishtan2011333}. It can be achieved with modifying both maps with Laplacian filter before computing the CCF. Motivation of its usage and the filter kernel are described in \cite{Chacon2002375}.

The ideas of CCF and LAP allows us to introduce four overlapping scores. 
Firstly, CCF itself can be computed. 
Then overlapping shows how incompatible are the positions of two proteins. 
For example, if they are too close or even has two atoms in a same position, the CCF will be big. 
This approach is denoted as CCF, and the goal is to minimize it.

Secondly, both maps can be filtered with the Laplacian filter to find the best match of their contours. 
This score is called the Contact score, and the goal is to maximize it.

The last two approaches imply applying the Laplacian filter to only one map, hence these scores shows how the contour of one map fits the volume of the other. 
These scores are called Skin-Core and Core-Skin scores and must be maximized to find the best match.

All scores except CCF are computed and then multiplied by $-1$ to write the optimization problem as the minimization one for all four cases. 
These scores allows us to introduce a matrix where each element considers overlapping between two positions of proteins.

\begin{Def} 
Let $\bQ \in \bbR^{n \times n}$ be an overlapping matrix that corresponds to a density map, where $n=m\cdot N$.
\end{Def}
Each matrix element $q_{ab}$ shows overlapping between $i$-th and $j$-th components of the complex, which are in $k$-th and $l$-th positions respectively, so $a=\left(i - 1\right)\cdot N + k$ and $b=\left( j - 1\right)\cdot N + l$. 
Overlapping between two positions of a single protein is set to zero.

Besides relative positions of proteins, their fitting to the density map should be considered. 
It means that it is more important to fit a protein's map to the complex' map contour rather than its volume, because we need to achieve the original position within the complex. 

\begin{Def}
Relevance of a protein's position is a measure of quality of it's fit into a density map's contour.
\end{Def}
Relevance can be measured with the LAP, since it can be used for contour matching.

With that, a relevance vector, which describes each possible position, can be introduced.

\begin{Def}
Let $\bb \in \bbR^n$ be a component relevance vector. Each element $b_{a}$ shows relevance of $i$-th component to $k$-th position in the complex, where $a=\left(i - 1\right)\cdot N + k$.
\end{Def}

Further, the problem's variable shows taken positions for each protein and is defined as following.
\begin{Def}
Let $\bx \in \{ 0, 1 \}^n$ be a binary vector that represents the proteins positions in the complex:
\[
x^k_i = 
\begin{cases}
1, &\mbox{if $i$-th protein is in the $k$-th position},\\
0 , &\mbox{otherwise},
\end{cases}
\]
where $a=\left(i - 1\right)\cdot N + k$.
\end{Def}
The vector $\bx$ is divided into $m$ subvectors for the proteins, each with length of $N$ for possible positions.

Since the best set of positions should have minimum overlapping between proteins and maximum relevance between each protein and the map, formulate the problem as a constrained binary quadratic optimization problem. It considers the overlapping in the quadratic term and the relevance in the linear term:
\begin{equation}
\begin{split}
\bx^* =&\argmin\limits_{\bx \in \{0,1\}^n} \left(\bx{\T}\bQ \bx - \bb{\T} \bx\right),\\
\text{s.t. }&\bA \bx = \mathbf{1}_m,
\end{split}
\label{eq::optim_problem}
\end{equation}
where $\mathbf{1}_m$ is a vector of ones and $\bA \in \bbR^{m \times Nm}$ is a matrix ensures that each protein within the complex takes a single position. Hence, it has the following structure:
\begin{equation}
\mathbf{A} = \left[\begin{array}{rrrrrrrrrr}
1 & \cdots & 1 & 0 & \cdots & 0 & \cdots & 0 & \cdots & 0 \\
0 & \cdots & 0 & 1 & \cdots & 1 & \cdots & 0 & \cdots & 0 \\
\vdots & \ddots & \vdots & \vdots & \ddots & \vdots & \cdots & \vdots & \ddots & \vdots \\
\undermat{N}{0 & \cdots & 0} & \undermat{N}{0 & \cdots & 0} & \cdots & \undermat{N}{1 & \cdots & 1} \\
\end{array}\right].
\end{equation}
$ $

Since matrix $\bQ$ is indefinite in general case, the problem~\eqref{eq::optim_problem} is non-convex. Moreover, the problem~\eqref{eq::optim_problem} is NP-hard. Because of that, relax integer constraints of the problem into continuous variables. The problem then can be solved with continuous optimization methods. The last approach will be discussed later, but reformulation in a continuous form can be written now as
\begin{equation}
\begin{split}
\by^* =&\argmin\limits_{\by \in [0,1]^n} \left(\by{\T}\bQ \by - \bb{\T} \by\right),\\
\text{s.t. }&\bA \by = \mathbf{1}_m.
\end{split}
\label{eq::optim_problem_contin}
\end{equation}

To define which proteins takes which place, return to the binary vector. 
To do that, replace the biggest element in each of $N$ position subvectors by 1 and others by 0:
\[
	i=1,\dots ,m:\ x_a = 
\begin{cases}
1, &\mbox{if }k=\argmax\limits_{k=1,\dots , N} y_a,\\
0 , &\mbox{otherwise}
\end{cases}
\]
for $a=\left(i - 1\right)\cdot N + k$.
\subsection{Convex relaxations}

Since \eqref{eq::optim_problem_contin} is a non-convex problem over a convex set, it can not be solved directly with guarantees of global minimum. However, an approximate solution can be found with relaxing the problem into a convex one. In the next sections spectrum shift and semidefinite relaxations of the problem are introduced. These approaches help in finding an approximate solution, using solutions of the relaxed problems.

\paragraph{Spectrum shift relaxation (Shift).} 
We shift the spectrum of matrix $\bQ$ to achieve positive-semidefiniteness and hence make the problem convex. 
The corresponding transformation is
\begin{equation}
\hat{\bQ} = \bQ - \lambda_{\min} \mathbf{I},
\label{eq::shift}
\end{equation}
where $\lambda_{\min}$ is the smallest eigenvalue of $\bQ$ and $\mathbf{I}$ is an identity matrix of the size of $\bQ$. 
Then the problem \eqref{eq::optim_problem_contin} can be rewritten with the new matrix $\hat{\bQ}$ as
\begin{equation}
\begin{split}
\by^* =&\argmin\limits_{\by \in [0,1]^n} \left(\by{\T}\hat{\bQ}\by - \bb{\T} \by\right),\\
\text{s.t. }&\bA \by = \mathbf{1}_m.
\end{split}
\label{eq::optim_problem_shift}
\end{equation}
The problem \eqref{eq::optim_problem_shift} can now be easily solved as a convex on. However, this method does not guarantee the global minimum of the initial problem~\eqref{eq::optim_problem}.


\paragraph{Semidefinite relaxation (SDP).} 
To introduce the semidefinite relaxation, rewrite the problem \eqref{eq::optim_problem_contin} as 
\begin{equation}
\begin{split}
\by^* =&\argmin\limits_{\by \in [0,1]^n} \left(\mathrm{Tr}\left(\bQ\bY\right) - \bb{\T} \by\right),\\
\text{s.t. }&\bA \by = \mathbf{1}_m,\\
&\bY=\by\by{\T}.
\end{split}
\end{equation}

To get a lower bond of the solution, relax the last constraint from equalities to inequalities, so now it is positive semidefinite: 
\[
\bY-\by\by{\T}\succeq 0.
\] 
But the initial binary program implies 
\[
\mathrm{diag}\left( \bx\bx{\T}\right)=\bx.
\] 
Hence, we use an additional constraint 
\[
\mathrm{diag}\left( \bY\right)=\by
\] 
to bound the problem. The relaxed problem is now convex and can be written as
\begin{equation}
\begin{split}
\by^* =&\argmin\limits_{\by \in [0,1]^n} \left(\mathrm{Tr}\left(\bQ\bY\right) - \bb{\T} \by\right),\\
\text{s.t. }&\bA \by = \mathbf{1}_m,\\
&\bY-\by\by{\T}\succeq 0,\\
&\mathrm{diag}\left( \bY\right)=\by.
\end{split}
\label{eq::optim_problem_sdp}
\end{equation}

\paragraph{Sequential quadratic programming (SQP).}
The basic ideas of sequential quadratic programming are described in Chapter 18 of the book \cite{nocedal2006numerical}. This approach finds a local minimum for a non-convex problem and implies solving a quadratic subproblem at each iteration. The subproblem is a convex second-order approximation of the Lagrangian function of the \eqref{eq::optim_problem_contin}, i.e. it involves a positive-semidefinite approximation of the Hessian. We use the implementation of an SQP algorithm from MATLAB Optimization Toolbox \cite{matlab2015b}. The initial point for the algorithm is obtained with solving the linear part of the problem \eqref{eq::optim_problem_contin}:
\begin{equation}
\begin{split}
\by^* =&\argmin\limits_{\by \in [0,1]^n} \left( - \bb{\T} \by\right),\\
\text{s.t. }&\bA \by = \mathbf{1}_m.
\end{split}
\label{eq::optim_problem_lin}
\end{equation}
In this work, the solution of the problem \eqref{eq::optim_problem_lin} can be treated as an approximation that only fits the density map, but does not consider overlapping between proteins.

\paragraph{Simulated annealing (SA).} The simulated annealing method \cite{Ingber199632} is a probabilistic global optimization method, implemented in MATLAB Global Optimization Toolbox \cite{matlab2015b}. This approach simulates a physical process of heating and then slow lowering the temperature of a material to decrease defects. At each iteration, it generates a new point near the current one, with a uniformly random direction and step length equals the current temperature. If the new point is better than the current one, the algorithm accepts it. If not, it accepts the point with probability
\begin{equation}
 \mathbb{P}\left(\mathrm{accept}\ x_{k+1}\right) =\left(1+\exp\left(\frac{\Delta}{T_k} \right)\right)^{-1},
\end{equation}
where $\Delta=f\left( x_{k+1}\right) -f\left( x_k\right)$ for an objective function $f$, current and new points $x_k$ and $x_{k+1}$ respectively and the current temperature $T_{k}$, which changes as
\[
T_{k+1}=0.95 \cdot T_k.
\]
For this work, the method can not be directly implemented for the problem \eqref{eq::optim_problem_contin} because it has linear constrains, but the method is designed for unconstrained  and bound-constrained problems. However, constraints can be implemented as a penalty function to the objective one, so the problem is converted as
\begin{equation}
\begin{split}
\by^* =&\argmin\limits_{\by \in [0,1]^n} \left(\by{\T}\bQ \by - \bb{\T} \by + w \| \bA \by - \mathbf{1}_m \|_1\right),\\
\end{split}
\label{eq::optim_problem_penalty}
\end{equation}
where $w\in\mathbb{R}$ is a penalty weight and $\| g\|_1$ denotes the $l_1$-norm of a vector $g$. The initial point for the algorithm is obtained by solving \eqref{eq::optim_problem_lin}. Moreover, since the algorithm has two parameters, the initial temperature $T_0$ and the penalty weight $w$, the method is denoted as SA$\left( T_0, w\right)$.
\subsection{Scoring functions to measure quality of fit} 

As the simplest scoring function which implies knowledge of the real structure of the protein, one can use root-mean-square deviation (RMSD) \cite{Maiorov1994625}. It measures the distance $\delta_i$ between pairs of $i$-th atoms of a protein in two positions, one in predicted position and one in the native position. Both atoms in a pair take the same place in a corresponding protein. With $M$ pairs of such atoms, it can be written as
\begin{equation}
\mathrm{RMSD}=\sqrt{\frac{1}{M}\sum_{i=1}^M\delta_i^2}.
\end{equation}
This approach helps to measure quality of other criteria on test data, but can not be used with direct determination of an unknown structure. 

Another approach is quality criteria that use only information about the density map. For future work, we propose two scoring functions recommended in the paper \cite{Vasishtan2011333}. Both of them use electron density maps of an initial structure and one from the solution. A way to produce a probe density map from the discrete solution $\bx^*$ is described in \cite{Vasishtan2011333}. It includes following steps.

\begin{enumerate}
\item[1.] Get the atomic structure from fitted proteins using the discrete solution $\bx^*$.
\item[2.] Impose a 3D grid with voxel size of 1 \AA .
\item[3.] For every non-hydrogen atom increase the density value of the nearest voxel by the atomic number of the atom.
\item[4.] Apply the Gaussian Fourier filter to blur the map. The recommended in \cite{Vasishtan2011333} sigma is $0.187\times$ resolution. The Gaussian kernel size is $2\cdot \lceil 2\sigma\rceil+1$, where $\lceil x\rceil$ is the smallest integer greater than or equal to $x$.
\item[5.] Resample the grid using Fourier method to match the sampling of the target map.
\end{enumerate}
When the probe map is obtained, compare the original (target) map and the probe one with the scoring functions describing above.

For $10$\AA$\ $resolution or less the authors propose the Laplacian-filtered cross-correlation function. The cross-correlation function itself is described above \eqref{eq::ccf}. For the LAP, modify both target and probe maps with Laplacian filter before computing the CCF.
For other cases, the mutual information score (MI) is proposed. The scoring function is
\begin{equation}
I\left( X, Y\right)=\underset{x\in X}{\sum}\underset{y\in Y}{\sum} p(x,y)\log{\frac{p(x,y)}{p(x)p(y)}}.
\label{eq::MI}
\end{equation}

Here, $X$ and $Y$ correspond to the density values in the probe and target maps. Functions $p(x)$ and $p(y)$ are the percentage of values in maps equal to $x$ and $y$. The aligned maps are maps where elements with equal coordinates represent one point in space. For aligned target and probe maps, $p(x,y)$ is percentage of elements with value $x$ in the probe map and $y$ in the target one. Because of wide range of values and noise, $X$ and $Y$ have limited number of values, e.g. 20 in \cite{Vasishtan2011333}.

\section{Dataset}

All methods were tested on simulated maps, since it allows one to compare the methods for different protein complexes that are tested in the same conditions, i.e. the same resolution of the map. The maps were generated as described below from 7 protein complexes. Their PDB entries are \textit{1e6v} \cite{1e6v}, \textit{1gte} \cite{1gte}, \textit{1tyq} \cite{1tyq}, \textit{1z5s} \cite{1z5s}, \textit{2p4n} \cite{2p4n}, \textit{4a6j} \cite{4a6j}, \textit{4bij} \cite{4bij}. The map resolution is 10 \AA , a voxel size is 1 \AA . Map's generation is similar to creating a map form a solution and includes following steps.
\begin{enumerate}
\item[1.] Impose a 3D grid with voxel size of 1 \AA .
\item[2.] For every non-hydrogen atom increase the density value of the nearest voxel by the atomic number of the atom.
\item[3.] Apply the Gaussian Fourier filter to blur the map. The recommended in \cite{Vasishtan2011333} sigma is $0.187\times$resolution. The Gaussian kernel size is $2\cdot \lceil 2\sigma\rceil+1$.
\end{enumerate}
\begin{table}

\end{table}

\section{Computational experiment}
Quality of fit is characterized by achieved RMSD of each protein in a complex~\eqref{tab::all_rmsd}. The solution is treated as correct if each protein has $\mathrm{RMSD}\leq 10$ \AA . Every results is characterized by correct answers ratio
\begin{equation*}
	\beta = \frac{N_c}{N},
\end{equation*}
where $N_c$ is a number of correctly determined protein's positions and $N$ is a number of proteins within the complex. Results for simulated annealing method are presented for following parameters: initial temperature is 100, penalty weight is 1. Results for other parameters were the same for all simulations (data not shown), which is connected with the quality of the linear solution described below. Moreover, only RMSD is used as a scoring function, since the purpose of the current work is to test the proposed optimization approach, and RMSD completely represents how precise an obtained solution is.
 
\paragraph{Comparison of quality of fit and the objective function's value.}
Results presented in Table \eqref{tab::1gte} for the complex \textit{1gte} show that the highest $\beta$ corresponds to the lowest approximate objective value. Moreover, the SQP and SA methods find a continuous solution, which gives almost the same objective value as the binary solution.
 
\paragraph{Comparison of the methods with the linear approximation.} Table \eqref{tab::all_rmsd} with results for \textit{1e6v}, \textit{1gte}, \textit{1z5s} shows that the highest $\beta$ was obtained by simulated annealing and sequential programming approaches. However, the initial point for both methods was obtained from the solution of the linear problem \eqref{eq::optim_problem_lin}, and $\beta$ for SQP, SA and the linear problem is almost the same. Hence, correct solutions for these datasets can be achieved using only the linear approach.

\paragraph{Overall performance.} Despite for some complexes mentioned above the linear approach leads to the correct solution (and to the correct solution with SQP and SA), in general it is not true. Moreover, SQP and SA sometimes make the linear solution even worse, which can be observed on the following Figure \ref{fig::overall}, that shows average $\beta$ for all complexes and map scores. Therefore, in these cases fitting to the given map is more important than arranging proteins' positions among each other.

\begin{table*}[p]
\centering
\caption{Objective function's optimal values for \textit{1gte}, Contact, CCF, Skin-Core and Core-Skin scores.}
\begin{tabular}{|c|c|c|c|c|c|c|c|}
\hline
&\multicolumn{3}{|c|}{Contact}&\multicolumn{3}{|c|}{CCF}\\
\hline
Method&Continuous&Binary&$\beta$&Continuous&Binary&$\beta$\\
\hline
SDP&-62711135&23648296& 0.50 &-2261626&563438& 0.25\\
\hline
Shift&1990252&1483316& 0.00 & 78968&43276& 0.00\\
\hline
SQP&-3707363&-3707363&1.00 &1086&1086& 0.00\\
\hline
SA$\left( 100, 1\right)$&-3707354&-3707363& 1.00&205462&205462&1.00\\
\hline
Linear&-80&-3707363& 1.00&-80&205462&1.00\\
\hline
&\multicolumn{3}{|c|}{Skin-Core}&\multicolumn{3}{|c|}{Core-Skin}\\
\hline
Method&Continuous&Binary&$\beta$&Continuous&Binary&$\beta$\\
\hline
SDP&-26429514&9343938&0.50&-26196482&2055248&0.25\\
\hline
Shift&801317&499195&0.00 &789572&210241&0.00\\
\hline
SQP&-2382769&-2382769&1.00&-2352385&-2352385& 1.00\\
\hline
SA$\left( 100,1\right)$&-2382764&-2382769&1.00&-2352380&-2352385& 1.00\\
\hline
Linear&-80&-2382769&1.00&-80&-2352385&1.00\\
\hline
\end{tabular}

\label{tab::1gte}
\end{table*}

\begin{table*}[p]
\centering
\caption{Results for three complexes. Each column represents a method, each value is a correct answers ratio $\beta$.} 
\begin{tabular}{|c|c|c|c|c|c|}
\hline
\multicolumn{6}{|c|}{\emph{1e6v}}\\
\hline
Scoring&SDP&Shift&SQP&SA&Linear\\
\hline
Contact&0.33&0.17&0.67&0.67&0.67\\
\hline
CCF&0.33&0.00&0.33&0.67&0.67\\
\hline
Skin-Core&0.50&0.00&0.67&0.67&0.67\\
\hline
Core-Skin&0.67&0.17&0.67&0.67&0.67\\
\hline
\multicolumn{6}{|c|}{\emph{1gte}}\\
\hline
Scoring&SDP&Shift&SQP&SA&Linear\\
\hline
Contact&0.50&0.00&1.00&1.00&1.00\\
\hline
CCF&0.25&0.00&0.00&1.00&1.00\\
\hline
Skin-Core&0.50&0.00&1.00&1.00&1.00\\
\hline
Core-Skin&0.25&0.00&1.00&1.00&1.00\\
\hline
\multicolumn{6}{|c|}{\emph{1z5s}}\\
\hline
Scoring&SDP&Shift&SQP&SA&Linear\\
\hline
Contact&0.25&0.00&0.50&1.00&1.00\\
\hline
CCF&0.25&0.00&0.25&1.00&1.00\\
\hline
Skin-Core&0.25&0.00&0.75&1.00&1.00\\
\hline
Core-Skin&0.25&0.00&0.75&1.00&1.00\\
\hline
\end{tabular}

\label{tab::all_rmsd}
\end{table*}

\begin{figure}
	\includegraphics[scale=0.45]{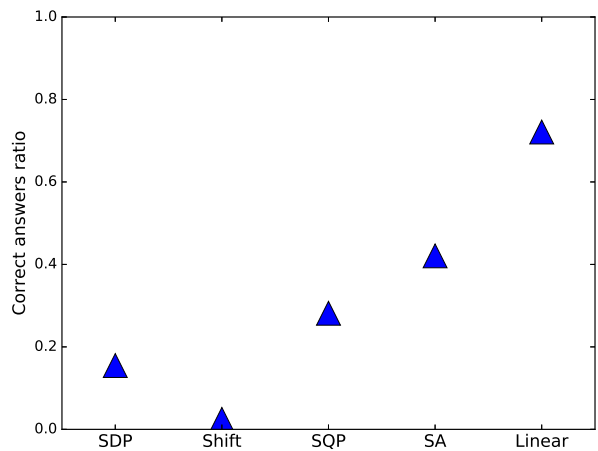}
	\caption{Averaged for all complexes and all map scores~$\beta$ }
	\label{fig::overall}
\end{figure}

\section{Conclusion}

The paper investigates quadratic programming approach for the problem of determining proteins' position inside a complex by its EM density map. The mathematical optimization problem is formulated using overlapping of proteins in computed positions between each other, which forms an overlapping matrix, and with the given math, which gives a relevance vector. 

We tested semidefinite relaxation, spectrum shift relaxation, sequential quadratic programming as quadratic programming approaches and simulated annealing and linear approximation for comparison with quadratic methods. The datasets were formed from simulated density maps. The best performance was shown by SQP, SA and linear approximation methods, which shows the importance of fitting a protein to the given map rather than looking for the best positions with relate to other proteins in the complex. 

\section{Acknowledgements}

This work is supported by RFBR, grant 16-37-00485.

\end{document}